# Accuracy Directly Controlled Fast Direct Solutions of General $\mathcal{H}^2$-Matrices and Its Application to Electrically Large Integral-Equation-Based Electromagnetic Analysis

Miaomiao Ma, *Student Member, IEEE,* and Dan Jiao, *Fellow, IEEE*

*Abstract*—The dense matrix resulting from an integral equation (IE) based solution of Maxwell's equations can be compactly represented by an $\mathcal{H}^2$-matrix. Given a general dense $\mathcal{H}^2$-matrix, prevailing fast direct solutions involve approximations whose accuracy can only be indirectly controlled. In this work, we propose new accuracy-controlled direct solution algorithms, including both factorization and inversion, for solving general $\mathcal{H}^2$-matrices, which does not exist prior to this work. Different from existing direct solutions, where the cluster bases are kept unchanged in the solution procedure thus lacking explicit accuracy control, the proposed new algorithms update the cluster bases and their rank level by level based on prescribed accuracy, without increasing computational complexity. Zeros are also introduced level by level such that the size of the matrix blocks computed at each tree level is the rank at that level, and hence being small. The proposed new direct solution has been applied to solve electrically large volume IEs whose rank linearly grows with electric size. A complexity of $O(NlogN)$ in factorization and inversion time, and a complexity of $O(N)$ in storage and solution time are both theoretically proven and numerically demonstrated. For constant-rank cases, the proposed direct solution has a strict $O(N)$ complexity in both time and memory. Rapid direct solutions of millions of unknowns can be obtained on a single CPU core with directly controlled accuracy.

*Index Terms*—Fast direct solvers, $\mathcal{H}^2$-matrix, dense matrices, fast solvers, factorization, inversion, electromagnetic analysis, frequency domain, electrodynamic, linear complexity solvers, integral equations, volume integral equations, electrically large analysis, three dimensional structures

## I. INTRODUCTION

THE $\mathcal{H}^2$-matrix is a general mathematical framework [1], [2], [3], [4] for compact representation and efficient computation of dense matrices, which can be utilized to develop fast solvers for electromagnetic analysis. Many existing fast solvers can be interpreted in this framework, although their developments may predate this framework. For example, the matrix structure resulting from the well-known FMM-based methods [5], [6] is an $\mathcal{H}^2$-matrix. In other words, the $\mathcal{H}^2$-matrix can be viewed as an algebraic generalization of the FMM method. Multiplying an $\mathcal{H}^2$-matrix by a vector has a complexity of $O(NlogN)$ for solving electrically large surface integral equations (SIEs). Most of the fast solvers accelerate the dense matrix-vector multiplication in an integral-equation based solution of electromagnetic problems, thus being fast iterative solvers. In mathematical literature like [7], it is shown that the storage, matrix-vector multiplication (MVM), and matrix-matrix multiplication of an $\mathcal{H}^2$-matrix can be performed in optimal $O(N)$ complexity for constant-rank cases. However, no such a complexity is shown for either matrix factorization or inversion, regardless of whether the rank is a bounded constant or an increasing variable. Later in [8], [9], fast matrix inversion and LU factorization algorithms are developed for $\mathcal{H}^2$-matrices. They have shown a complexity of $O(N)$ for solving constant-rank cases like electrically small and moderate problems for both surface and volume IEs (VIEs) [10], [11], [12], [13]; and a complexity of $O(NlogN)$ for solving electrically large VIEs [14], [15].

Despite a significantly reduced complexity, in existing direct solutions of $\mathcal{H}^2$-matrices like [8], [9], the cluster bases used to represent a dense matrix are also used to represent the matrix's inverse as well as LU factors. During the inversion and factorization procedure, only the coupling matrix of each admissible block needs to be computed, while the cluster bases are kept to be the same as those in the original matrix. Physically speaking, such a choice of the cluster bases for the inverse matrix often constitutes an accurate choice. However, mathematically speaking, the accuracy of the inversion and factorization cannot be directly controlled. The lack of a direct accuracy control is also observed in the $\mathcal{H}^2$-based matrix-matrix multiplication in mathematical literature such as [7]. When the accuracy of the direct solution is not satisfactory, one can only change the original $\mathcal{H}^2$-representation, i.e., change the cluster bases or the rank for representing the original matrix, with the hope that they can better represent the inverse and LU factors. Therefore, the accuracy of the resulting direct solution is not directly controlled.

Recently, an HSS-matrix structure [16] has also been explored to develop fast direct solutions of exact arithmetic [17], [18] for electromagnetic analysis. In other words, the direct solution of the HSS matrix does not involve any approximation. However, the HSS matrix, as a special class of $\mathcal{H}^2$-matrix, requires only one admissible block be formed for each node in the $\mathcal{H}^2$-tree. All the off-diagonal matrices in the HSS matrix are represented as low-rank matrices. As a result, the resultant rank can be very high to be computed efficiently, especially for electrically large analysis. This is because as

Manuscript received on Mar. 18, 2017. This work was supported by a grant from NSF under award No. 1619062, a grant from SRC (Task 1292.073), and a grant from DARPA under award HR0011-14-1-0057. The authors are with the School of Electrical and Computer Engineering, Purdue University, West Lafayette, IN, 47907 USA (e-mail: djiao@purdue.edu).

long as the sources and observers are separated, even though their distance is infinitely small, their interaction has to be represented by a low-rank matrix to suite the structure of an HSS matrix. From this perspective, an $\mathcal{H}^2$-matrix is a more efficient structure for general applications, since the distance between the sources and observers can be used to reduce the rank required to represent a dense matrix for a prescribed accuracy. However, the accuracy-controlled direct solution of general $\mathcal{H}^2$-matrices is still lacking in the open literature. The contribution of this work is such an accuracy-controlled direct solution of general $\mathcal{H}^2$-matrices. This solution has a directly controlled accuracy, while preserving the low computational complexity of the $\mathcal{H}^2$-matrices. This solution can be applied to solve a general $\mathcal{H}^2$-matrix, be it real- or complex-valued, symmetric or unsymmetrical. As a demonstration, we show how it is used to solve electrically large VIEs, with an $O(NlogN)$ complexity in factorization and inversion, and $O(N)$ complexity in solution and memory, and with a strictly controlled accuracy.

The rest of the paper is organized as follows. In Section II, we review the mathematical preliminaries. In Section III, we present the proposed factorization and inversion algorithms for general $\mathcal{H}^2$-matrices. In Section IV, we describe the proposed matrix solution algorithm including both backward and forward substitutions. The accuracy and computational complexity of the proposed direct solutions are analyzed in Section V. Section VI demonstrates the application of this work in solving electrically large VIEs. In Section VII, we summarize this paper.

## II. PRELIMINARIES

An $\mathcal{H}^2$-matrix [1], [2], [3], [4] is generally stored in a tree structure, with each node in the tree called a cluster, as shown by the left or right tree illustrated in Fig. 1a. This tree can be obtained by recursively dividing the entire unknown set into two subsets until the number of unknowns in each subset is no greater than $leafsize$, a predefined constant. The tree built for row unknowns of a matrix can be termed as a row tree, while that for column unknowns is called a column tree. In a Galerkin-based scheme for obtaining a system matrix, the two trees are the same. Using the row- and column tree, the original matrix can be partitioned into multilevel admissible and inadmissible blocks. This is done by checking the following admissibility condition level by level between a cluster $t$ in the row tree and a cluster $s$ in the column tree

$$max\{diam(\Omega_t), diam(\Omega_s)\} \leq \eta dist(\Omega_t, \Omega_s), \quad (1)$$

where $\Omega_t$ ($\Omega_s$) denotes the geometrical support of the unknown set $t$ ($s$), $diam\{\cdot\}$ is the Euclidean diameter of a set, $dist\{\cdot,\cdot\}$ denotes the Euclidean distance between two sets, and $\eta$ is a positive parameter that can be used to control the admissibility condition. If the two clusters are admissible, i.e., satisfying the above condition, their children clusters do not need to be checked. As an example, a four-level block cluster $\mathcal{H}^2$-tree is illustrated in Fig. 1a, where a green link denotes an admissible block, formed between a row cluster and a column cluster that satisfy (1); and a red link denotes an inadmissible block. The resultant $\mathcal{H}^2$-matrix is shown in Fig. 1b.

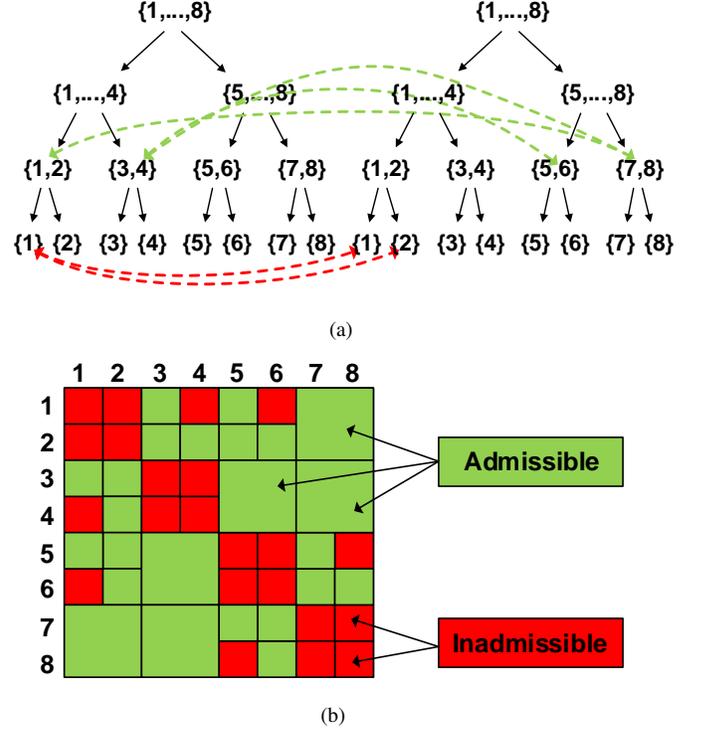

Fig. 1: Illustration of a block cluster tree and resulting $\mathcal{H}^2$-matrix partition. (a) Block cluster tree. (b) $\mathcal{H}^2$-matrix structure.

An admissible block in an $\mathcal{H}^2$-matrix is represented as

$$\mathbf{Z}_{t,s} = (\mathbf{V}_t)_{\#t \times k} (\mathbf{S}_{t,s})_{k \times k} (\mathbf{V}_s)^T_{\#s \times k}, \quad (2)$$

where $\mathbf{V}_t$ ($\mathbf{V}_s$) is called cluster basis of cluster $t$ ($s$), $\mathbf{S}_{t,s}$ is called coupling matrix. The cluster basis in an $\mathcal{H}^2$-matrix has a nested property. This means the cluster basis for a non-leaf cluster $t$, $\mathbf{V}_t$, can be expressed by its two children's cluster bases, $\mathbf{V}_{t_1}$ and $\mathbf{V}_{t_2}$, as

$$(\mathbf{V}_t)_{\#t \times k} = \begin{bmatrix} (\mathbf{V}_{t_1})_{\#t_1 \times k_1} & 0 \\ 0 & (\mathbf{V}_{t_2})_{\#t_2 \times k_2} \end{bmatrix} \begin{bmatrix} (\mathbf{T}_{t_1})_{k_1 \times k} \\ (\mathbf{T}_{t_2})_{k_2 \times k} \end{bmatrix} \quad (3)$$

where $\mathbf{T}_{t_1}$ and $\mathbf{T}_{t_2}$ are called transfer matrices. In an $\mathcal{H}^2$-matrix, the number of blocks formed by a single cluster at each tree level is bounded by a constant, $C_{sp}$. The size of an inadmissible block is $leafsize$, and hence inadmissible blocks appear only at the leaf level. Because of the nested relationship of cluster bases, the cluster bases only need to be stored for leaf clusters. For non-leaf clusters, only transfer matrices need to be stored.

In mathematical literature [7], for constant-rank cases, it is shown that the computational complexity of an $\mathcal{H}^2$-matrix is optimal (linear) in storage, matrix-vector multiplication, and matrix-matrix multiplication. In [8], [9], [10], [11], [12], [13], it is shown that $O(N)$ direct solutions can also be developed for $\mathcal{H}^2$-matrices for constant-rank cases as well as electrically moderate cases whose rank can be controlled by a rank function. For electrically large analysis, the rank of an $\mathcal{H}^2$-representation of IE kernels is not constant any more for achieving a prescribed accuracy. Instead, the rank for an off-diagonal interaction is shown to scale linearly with the electric size for general 3-D configurations [21]. Such a rank's growth rate with electric size should be taken into account in developing low-complexity direct solutions for electrically large analysis. Recent work on direct solutions of the FMM can also be seen in [19], [20].

## III. PROPOSED FACTORIZATION AND INVERSION ALGORITHMS

Given a general $\mathcal{H}^2$-matrix $\mathbf{Z}$, the proposed direct solution is a one-way tree traversal from the leaf level all the way up to the root level of the $\mathcal{H}^2$-tree. At each level, the factorization proceeds from the left to the right across all the clusters. This is very different from the algorithm of [9], [11], [14], [15], where a recursive procedure is performed. To help better understand the proposed algorithm, we will use the $\mathcal{H}^2$-matrix shown in Fig. 2(a) as an example to illustrate the overall procedure. In this figure, green blocks are admissible blocks, and red ones represent inadmissible blocks.

### A. Computation at the leaf level

*1) For the first leaf cluster:* We start from the leaf level ($l = L$). For the first cluster, whose index is denoted by $i = 1$, we perform three steps as follows.

**Step 1.** We compute the complementary cluster basis $\mathbf{V}_i^\perp$, which is orthogonal to $\mathbf{V}_i$ (cluster basis of cluster $i$). This can be easily done by carrying out an SVD or QL factorization on $\mathbf{V}_i$. The cost is not high, instead, it is constant at the leaf level. If $\mathbf{V}_i^\perp$ has already been generated and stored during the procedure of the $\mathcal{H}^2$-matrix construction, this step can be skipped. We then combine $\mathbf{V}_i$ with $\mathbf{V}_i^\perp$ to form $\widetilde{\mathbf{Q}}_i$ matrix as the following

$$\widetilde{\mathbf{Q}}_i = [(\mathbf{V}_i^\perp)_{\#i \times (\#i - k_l)}\ (\mathbf{V}_i)_{\#i \times k_l}], \quad (4)$$

in which $k_l$ is the rank at level $l$ (now $l = L$), and $\#i$ denotes the number of unknowns contained in cluster $i$.

**Step 2.** Let $\mathbf{Z}_{cur}$ be the current matrix that remains to be factorized. For leaf cluster $i = 1$, $\mathbf{Z}_{cur} = \mathbf{Z}$. We compute a transformed matrix $\mathbf{Q}_i^H \mathbf{Z}_{cur} \overline{\mathbf{Q}}_i$, and use it to overwrite $\mathbf{Z}_{cur}$, thus

$$\mathbf{Z}_{cur} = \mathbf{Q}_i^H \mathbf{Z}_{cur} \overline{\mathbf{Q}}_i, \quad (5)$$

where

$$\mathbf{Q}_i = diag\{\mathbf{I}_1, \mathbf{I}_2, ..., \widetilde{\mathbf{Q}}_i, ..., \mathbf{I}_{\#\{\text{leaf clusters}\}}\}, \quad (6)$$

is a block diagonal matrix whose $i$-th diagonal block is $\widetilde{\mathbf{Q}}_i$, and other blocks are identity matrices $\mathbf{I}$, the size of each of which is the corresponding leaf-cluster size. The subscripts in the right hand side of (6) denote the indexes of diagonal blocks, which are also the indexes of leaf clusters, and $\#\{\text{leaf clusters}\}$ stands for the number of leaf clusters. If leaf cluster $i = 1$ is being computed, $\widetilde{\mathbf{Q}}_i$ appears in the first block of (6). In addition, the $\overline{\mathbf{Q}}_i$ in (5) denotes the complex conjugate of $\mathbf{Q}_i$, as the cluster bases we construct are complex-valued to account for general $\mathcal{H}^2$-matrices. In (6), $\mathbf{Q}_i$ only consists of one $\widetilde{\mathbf{Q}}_i$, where $i$ is the current cluster being computed. This is because $\widetilde{\mathbf{Q}}_i$ for other clusters are not ready for use yet, since in the proposed algorithm, we will update cluster bases to take into account the contribution from fill-ins. This will become clear in the sequel.

The purpose of doing this step, i.e., obtaining (5), is to introduce zeros in the matrix being factorized, since

$$\widetilde{\mathbf{Q}}_i^H \times \mathbf{V}_i = \begin{bmatrix} \mathbf{0}_{(\#i - k_l) \times k_l} \\ \mathbf{I}_{k_l \times k_l} \end{bmatrix} \quad (7)$$

$$\mathbf{V}_i^T \times \overline{\widetilde{\mathbf{Q}}}_i = [\mathbf{0}_{k_l \times (\#i - k_l)}\ \mathbf{I}_{k_l \times k_l}]. \quad (8)$$

The operation of (5) thus will zero out the first $(\#i - k_l)$ rows of the admissible blocks formed by row cluster $i$, as well as the first $(\#i - k_l)$ columns in $\mathbf{Z}_{cur}$, as shown by the white blocks in Fig. 2(b). The real computation of (5) is, hence, in the inadmissible blocks of cluster $i$, as shown by the four red blocks associated with cluster $i = 1$ in Fig. 2(b).

**Step 3.** After Step 2, the first $(\#i - k_l)$ rows and columns of $\mathbf{Z}_{cur}$ in the admissible blocks of cluster $i$ become zero. Hence, if we eliminate the first $(\#i - k_l)$ unknowns of cluster $i$ via a partial LU factorization, the admissible blocks of cluster $i$ at these rows and columns would stay as zero, and thus not affected. As a result, the resulting L and U factors, as well as the fill-in blocks generated, would only involve a constant number of blocks, which is bounded by $O(C_{sp}^2)$.

In view of the above property, in Step 3, we perform a partial LU factorization to eliminate the first $(\#i - k_l)$ unknowns of cluster $i$ in $\mathbf{Z}_{cur}$. Let this set of unknowns be denoted by $i'$. $\mathbf{Z}_{cur}$ can be cast into the following form accordingly

$$\mathbf{Z}_{cur} = \begin{bmatrix} \mathbf{F}_{i'i'} & \mathbf{Z}_{i'c} \\ \mathbf{Z}_{ci'} & \mathbf{Z}_{cc} \end{bmatrix}, \quad (9)$$

in which $\mathbf{F}_{i'i'}$ denotes the diagonal block of the $(\#i - k_l)$ unknowns, $\mathbf{Z}_{i'c}$ is the remaining rectangular block in the row of set $i'$. Here, $c$ contains the rest of the unknowns which do not belong to $i'$. The $\mathbf{Z}_{ci'}$ is the transpose block of $\mathbf{Z}_{i'c}$, and $\mathbf{Z}_{cc}$ is the diagonal block formed by unknowns in set $c$. It is evident that $\mathbf{Z}_{i'c}$ ($\mathbf{Z}_{ci'}$) only consists of $O(C_{sp})$ inadmissible blocks, as the rest of the blocks are zero. Take leaf cluster $i = 1$ as an example. It forms inadmissible blocks with clusters $i = 2, 4, 6$, as shown by the red blocks in Fig. 2(b). Hence, $\mathbf{Z}_{i'c} = [\mathbf{F}_{12}, \mathbf{F}_{14}, \mathbf{F}_{16}]$.

After performing a partial LU factorization to eliminate the unknowns contained in $i'$, we obtain

$$\mathbf{Z}_{cur} = \begin{bmatrix} \mathbf{L}_{i'i'} & \mathbf{0} \\ \mathbf{Z}_{ci'}\mathbf{U}_{i'i'}^{-1} & \mathbf{I} \end{bmatrix} \begin{bmatrix} \mathbf{I} & \mathbf{0} \\ \mathbf{0} & \widetilde{\mathbf{Z}}_{cc} \end{bmatrix} \begin{bmatrix} \mathbf{U}_{i'i'} & \mathbf{L}_{i'i'}^{-1}\mathbf{Z}_{i'c} \\ \mathbf{0} & \mathbf{I} \end{bmatrix}, \quad (10)$$

where

$$\mathbf{F}_{i'i'} = \mathbf{L}_{i'i'}\mathbf{U}_{i'i'}, \quad (11)$$

and the Schur complement is

$$\widetilde{\mathbf{Z}}_{cc} = \mathbf{Z}_{cc} - \mathbf{Z}_{ci'}\mathbf{U}_{i'i'}^{-1}\mathbf{L}_{i'i'}^{-1}\mathbf{Z}_{i'c}, \quad (12)$$

in which $\mathbf{Z}_{ci'}\mathbf{U}_{i'i'}^{-1}\mathbf{L}_{i'i'}^{-1}\mathbf{Z}_{i'c}$ represents the fill-in blocks introduced due to the elimination of $i'$-unknowns. Because of the zeros introduced in the admissible blocks, the $\mathbf{Z}_{i'c}$ and $\mathbf{Z}_{ci'}$ only consist of inadmissible blocks formed by cluster $i$, the number of which is bounded by constant $C_{sp}$. Hence, the number of fill-in blocks introduced by this step is also bounded by a constant, which is $O(C_{sp}^2)$.

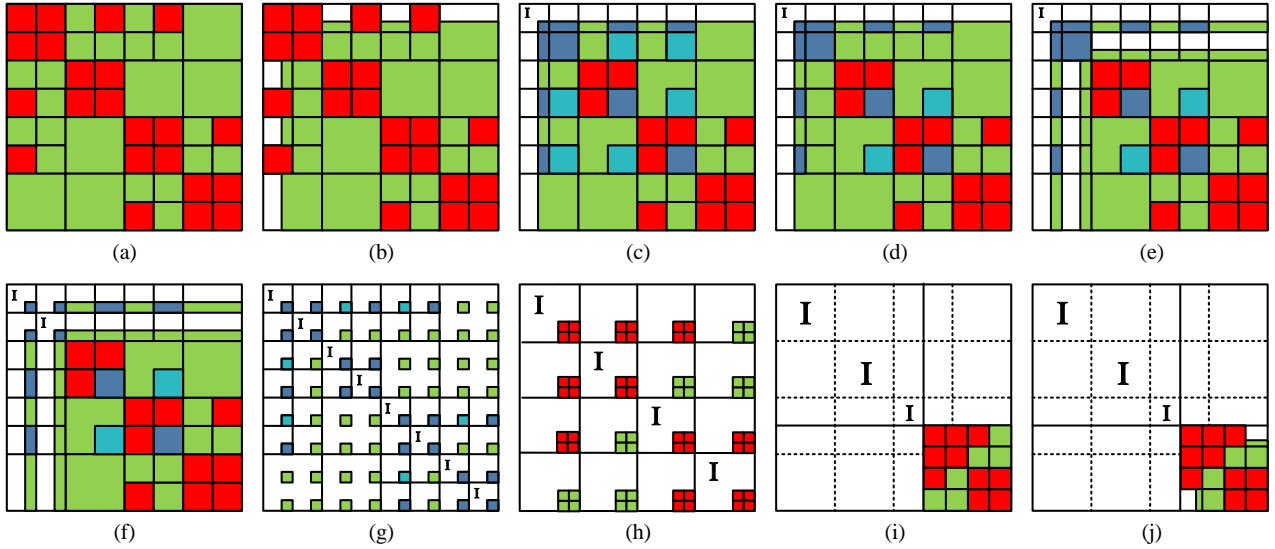

Fig. 2: Illustration of the overall procedure of the proposed direct solution. (a) An $\mathcal{H}^2$-matrix $\mathbf{Z}_{cur}$. (b) $\mathbf{Q}_1^H \mathbf{Z}_{cur} \overline{\mathbf{Q}}_1$. (c) Partial LU factorization for cluster 1. (d) Update cluster basis of cluster 2. (e) $\mathbf{Q}_2^H \mathbf{Z}_{cur} \overline{\mathbf{Q}}_2$. (f) Partial LU factorization for cluster 2. (g) $\mathcal{H}^2$-matrix after leaf level elimination. (h) Merging to next level. (i) Ancestor level structure. (j) Ancestor level $\mathbf{Q}_t^H \mathbf{Z}_{cur} \overline{\mathbf{Q}}_t$.

Eqn. (10) can further be written in short as

$$\mathbf{Z}_{cur} = \mathbf{L}_i^l \begin{bmatrix} \mathbf{I} & \mathbf{0} \\ \mathbf{0} & \widetilde{\mathbf{Z}}_{cc} \end{bmatrix} \mathbf{U}_i^l, \quad (13)$$

in which the first matrix of (10), which is the $\mathbf{L}$ factor generated after a partial LU of cluster $i$ at level $l$, is denoted by $\mathbf{L}_i^l$. Similarly, the rightmost matrix of (10), which is the $\mathbf{U}$ factor generated after a partial LU of cluster $i$ at level $l$, is denoted by $\mathbf{U}_i^l$. The center matrix is the one that remains to be factorized. Thus, we let

$$\mathbf{Z}_{cur} = \begin{bmatrix} \mathbf{I} & \mathbf{0} \\ \mathbf{0} & \widetilde{\mathbf{Z}}_{cc} \end{bmatrix}. \quad (14)$$

This matrix is illustrated in Fig. 2(c), where the fill-in blocks are added upon previous $\mathbf{Z}_{cur}$. If the fill-in block is added upon an inadmissible block, we add the fill-in directly to the inadmissible block. For the example shown in Fig. 2(a), these blocks are the blocks marked in dark blue in Fig. 2(c) after the $i'$ set of unknowns in cluster $i = 1$ is eliminated. If the fill-in appears in an admissible block, we store the fill-in block temporarily in this admissible block for future computation. For the example shown in Fig. 2(a), these blocks are the blocks marked in light blue in Fig. 2(c).

From (12), it can be seen that the fill-in block formed between cluster $j$ and cluster $k$ has the following form

$$\mathbf{F}_{j,k} = \mathbf{Z}_{ji'} \mathbf{U}_{i'i'}^{-1} \mathbf{L}_{i'i'}^{-1} \mathbf{Z}_{i'k}, \quad (15)$$

which is also low rank in nature. The row cluster $j$ and column cluster $k$ belong to the set of clusters that form inadmissible blocks with cluster $i$, whose number is bounded by $C_{sp}$. If the original block formed between $j$ and $k$, $\mathbf{Z}_{cur,jk}$, is admissible, no computation is needed for $\mathbf{F}_{j,k}$, i.e., we do not add it directly to the admissible block. Instead, we store it with this admissible block temporarily until cluster $j$ is reached during the factorization. At that time, $\mathbf{F}_{j,k}$ will be used to compute an updated cluster basis for $j$. If the elimination of different clusters results in more than one fill-in block in admissible $\mathbf{Z}_{cur,jk}$ block, the later ones are added upon the existing $\mathbf{F}_{j,k}$. In other words, the fill-ins at the admissible $\mathbf{Z}_{cur,jk}$ block are summed up and stored in a single $\mathbf{F}_{j,k}$.

In summary, Step 3 is to perform a partial LU factorization as shown in (10) to eliminate the first $\#i - k_l$ unknowns in cluster $i$. The storage and computation cost of this step can be analyzed as follows. In (10), the LU factorization of $\mathbf{F}_{i'i'}$ into $\mathbf{L}_{i'i'}$ and $\mathbf{U}_{i'i'}$ costs $O(leafsize^3)$ in computation, and $O(leafsize^2)$ in memory, which is hence constant. $\mathbf{Z}_{ci'}$ ($\mathbf{Z}_{i'c}$) only involves $O(C_{sp})$ inadmissible blocks. Hence, the computation of $\mathbf{Z}_{ci'} \mathbf{U}_{i'i'}^{-1}$ as well as $\mathbf{L}_{i'i'}^{-1} \mathbf{Z}_{i'c}$ is simply $O(C_{sp})$ operations, each of which scales as $O(leafsize^3)$, and hence the cost is constant also. The results of the partial factorization, i.e., $\mathbf{L}_i^l$ and $\mathbf{U}_i^l$, are stored in $O(C_{sp})$ blocks of $leafsize$. Thus, the memory cost is also constant. As for the fill-in blocks, their number is bounded by $C_{sp}^2$, and each of which is of $leafsize$. Hence, the computation and storage associated with fill-in blocks are constant too, at the leaf level for each cluster.

*2) For other leaf clusters:* Next, we proceed to the second leaf cluster $i = 2$ and other leaf clusters. For these leaf clusters, all the steps are the same as those performed for the first leaf cluster, except for adding another step before Step 1. We term this step as **Step 0**.

**Step 0.** This step is to update cluster basis $\mathbf{V}_i$ to account for the contributions of the fill-ins due to the elimination of previous clusters. If we do not update cluster bases $\mathbf{V}_i$, they are not accurate any more to represent the admissible blocks formed by cluster $i$, since the contents of these blocks have been changed due to the addition of the fill-in blocks. However, if we let the admissible block updated by fill-ins become a full-matrix block, i.e., an inadmissible block, the number of fill-ins would keep increasing after each step of partial factorization, and hence the resultant complexity would be too high to tolerate. To overcome this problem, we propose to update cluster bases as follows.

For cluster $i$, we take all the fill-in blocks associated with its admissible blocks, and use them to update cluster basis $\mathbf{V}_i$. These fill-in blocks are shown in (15), which have been generated and stored in the corresponding admissible blocks during the elimination of previous clusters. Let them be $\mathbf{F}_{i,j_k}(k = 1, 2, ..., O(C_{sp}))$, where the first subscript denotes the row cluster $i$, and the second being the column cluster index. The number of such blocks is bounded by $O(C_{sp})$ because the number of admissible blocks formed by a single cluster, like $i$, is no greater than $C_{sp}$. Next, we compute

$$\mathbf{G}_i = (\mathbf{I} - \mathbf{V}_i\mathbf{V}_i^H)(\sum_k \mathbf{F}_{i,j_k}\mathbf{F}_{i,j_k}^H)(\mathbf{I} - \mathbf{V}_i\mathbf{V}_i^H)^H. \quad (16)$$

We then perform an SVD on $\mathbf{G}_i$, and truncate the singular vector matrix based on prescribed accuracy $\epsilon_{fill-in}$, obtaining

$$\mathbf{G}_i \stackrel{\epsilon_{fill-in}}{=} (\mathbf{V}_i^{add})\Sigma(\mathbf{V}_i^{add})^H. \quad (17)$$

The $\mathbf{V}_i^{add}$ is the additional cluster basis we supplement, since it captures the part of the column space of the fill-in blocks, which is not present in the original $\mathbf{V}_i$. The cluster basis can then be updated for cluster $i$ as

$$\widetilde{\mathbf{V}}_i = [\mathbf{V}_i \ \mathbf{V}_i^{add}]. \quad (18)$$

By doing so, we are able to keep the original admissible blocks to be admissible while accurately taking into account the fill-in's contributions. Hence, the number of fill-ins introduced due to the elimination of each cluster can be kept to be a constant, instead of growing. In addition, by keeping the original cluster basis $\mathbf{V}_i$ in the new bases instead of generating a completely new one, we can keep the nested relationship in the original $\mathcal{H}^2$-matrix for efficient computations at upper levels. The computational cost of such a supplement of cluster bases is constant at leaf level, and scales as $O(k_l^3)$ at other levels (this will become clear soon) for each cluster.

In Fig. 2(d), the Step 0 is reflected by turning the light blue blocks associated with cluster $i = 2$ in Fig. 2(c) to green blocks. After $\mathbf{V}_i$ is updated to $\widetilde{\mathbf{V}}_i$, we find its complementary bases, i.e., perform Step 1, to obtain

$$\widetilde{\mathbf{Q}}_i = [(\widetilde{\mathbf{V}}_i^\perp)_{\#i \times (\#i - k_l)} \ (\widetilde{\mathbf{V}}_i)_{\#i \times k_l}]. \quad (19)$$

We then use $\widetilde{\mathbf{Q}}_i$ to generate (6), and project the current matrix to be factorized, shown in Fig. 2(d), to the matrix shown in (5). This again will zero out the first $(\#i - k_l)$ rows (columns) of the admissible blocks formed by cluster $i$, as illustrated in Fig. 2(e). We then proceed to Step 3 to perform a partial LU factorization. After this step, the matrix $\mathbf{Z}_{cur}$, which remains to be factorized, is illustrated in Fig. 2(f). Steps 0 to 3 are then repeated for the rest of the leaf clusters.

*3) Updating the coupling matrix at the leaf level and the transfer matrix at one level higher:* After the computation at the leaf level is finished, before proceeding to the next level, we need to update the coupling matrices at the leaf level, as well as the transfer matrices at one level higher. This is because the cluster bases have been updated at the leaf level.

For admissible blocks without fill-ins, their coupling matrices can be readily updated by appending zeros. For example, considering an admissible block formed between clusters $i$ and $k$. Its coupling matrix $\mathbf{S}_{i,k}$ is updated to

$$\widetilde{\mathbf{S}}_{i,k} = \begin{bmatrix} \mathbf{S}_{i,k} & \mathbf{0} \\ \mathbf{0} & \mathbf{0} \end{bmatrix}, \quad (20)$$

because in this way

$$\widetilde{\mathbf{V}}_i\widetilde{\mathbf{S}}_{i,k}\widetilde{\mathbf{V}}_k^T = \mathbf{V}_i\mathbf{S}_{i,k}\mathbf{V}_k^T. \quad (21)$$

For admissible blocks with fill-ins, their coupling matrices are updated by adding the fill-in part onto the augmented coupling matrix (shown in (20)) as

$$\widetilde{\mathbf{S}}_{i,k} = \widetilde{\mathbf{S}}_{i,k} + \widetilde{\mathbf{V}}_i^H \mathbf{F}_{i,k} \overline{\widetilde{\mathbf{V}}}_k, \quad (22)$$

in which $\mathbf{F}_{i,k}$ represents the fill-in block associated with the admissible block.

We also need to update the transfer matrices at one level higher, i.e., at the $l = L - 1$ level. This can be readily done by adding zeros to the original transfer matrices. To be clearer, consider a non-leaf cluster $t$, whose two children clusters are $t_1$ and $t_2$, its cluster basis can be written as

$$\begin{bmatrix} \mathbf{V}_{t_1}\mathbf{T}_{t_1} \\ \mathbf{V}_{t_2}\mathbf{T}_{t_2} \end{bmatrix} = \begin{bmatrix} [\mathbf{V}_{t_1}\mathbf{V}_{t_1}^{add}] & \mathbf{0} \\ \mathbf{0} & [\mathbf{V}_{t_2}\mathbf{V}_{t_2}^{add}] \end{bmatrix} \begin{bmatrix} \mathbf{T}_{t_1} \\ \mathbf{0} \\ \mathbf{T}_{t_2} \\ \mathbf{0} \end{bmatrix}, \quad (23)$$

where $\mathbf{T}_{t_1}$ and $\mathbf{T}_{t_2}$ are original transfer matrices associated with non-leaf cluster $t$.

The transfer matrices and coupling matrices at other levels all remain the same as before. They are not updated until corresponding levels are reached.

### B. Computation at the non-leaf level

After all the leaf-level computation, we obtain a matrix shown in Fig. 2(g), which is current $\mathbf{Z}_{cur}$ left to be factorized. All the empty blocks are zero blocks. The diagonal blocks corresponding to the unknowns that have been eliminated are identity matrices shown as $\mathbf{I}$. The nonzero blocks, shown as shaded blocks, each is of size $(k_{l+1} \times k_{l+1})$, where $l$ is current tree level, $l + 1$ is the previous tree level whose computation is finished, and $k_{l+1}$ denotes the rank at the $(l + 1)$-th level. The structure shown in Fig. 2(g) is not only true for $l = L - 1$ level, i.e., one level above the leaf level, but also true for every non-leaf level $l$, except that the size of each block is different at different tree levels, since $k_{l+1}$ is different especially for electrically large analysis. The computation at a non-leaf level is performed as follows.

*1) Merging and permuting:* In $\mathbf{Z}_{cur}$, whose structure is shown in Fig. 2(g), we permute the matrix, and merge the four blocks of a non-leaf cluster, formed between its two children clusters and themselves, to a single block, as shown by the transformation from Fig. 2(g) to Fig. 2(h). After this, we further permute the unknowns to group the unknowns left to be eliminated together, and order them at the end. This results in a matrix shown in Fig. 2(i). As a result, the bottom right block becomes the matrix that remains to be factorized. This matrix is of size $2^l(2k_{l+1} \times 2k_{l+1})$ by $2^l(2k_{l+1} \times 2k_{l+1})$, because it is formed between $2^l$ clusters at the non-leaf level $l$, and each block is of size $2k_{l+1} \times 2k_{l+1}$. Obviously, this matrix is again an $\mathcal{H}^2$-matrix, as shown by green admissible blocks and red inadmissible blocks in Fig. 2(i). The difference with the leaf level is that each block, be its inadmissible or admissible, now is of size $2k_{l+1} \times 2k_{l+1}$. Hence, in the proposed direct solution algorithm, at each non-leaf level, the computation is performed on the matrix blocks whose size is the *rank* at that level, and hence efficient.

Since each admissible block is only of size $2k_{l+1} \times 2k_{l+1}$, its expression is different from that at the leaf level. Consider an admissible block formed between clusters $t$ and $s$ at a non-leaf level $l$, as shown by the green blocks in Fig. 2(h). It has the following form:

$$(\text{Admissible block})^{(l)}_{t,s} = \mathbf{T}_t \mathbf{S}_{t,s} \mathbf{T}_s^T, \quad (24)$$

where $\mathbf{T}_t$ is the transfer matrix for non-leaf cluster $t$, whose size is $O(2k_{l+1} \times k_l)$, $\mathbf{T}_s$ is the transfer matrix for cluster $s$ whose size is also $O(2k_{l+1} \times k_l)$, and $\mathbf{S}_{t,s}$ is the coupling matrix whose size is $O(k_l \times k_l)$.

Since the matrix block left to be factorized after merging and permutation is again an $\mathcal{H}^2$-matrix, we can apply the same procedure performed at the leaf level, i.e., from Step 0 to Step 3, to the non-leaf level. In addition, from (24), it can also be seen clearly that now, transfer matrix $\mathbf{T}$ plays the same role as cluster basis $\mathbf{V}$ at the leaf level. Therefore, wherever $\mathbf{V}$ is used in the leaf-level computation, we replace it by $\mathbf{T}$ at a non-leaf level. In the following, we will briefly go through the Steps 0–3 performed at the non-leaf level, as the essential procedure is the same as that in the leaf level.

*2) Steps 0 to 3:* **Step 0 at the non-leaf level.** This step is to update transfer matrix $\mathbf{T}_i$ to account for the contributions of the fill-ins due to the elimination of previous clusters.

For cluster $i$, we take all the fill-in blocks associated with its admissible blocks, and use them to update transfer matrix $\mathbf{T}_i$. These fill-in blocks are shown in (15), which have been generated and stored during the factorization of previous clusters. Let them be $\mathbf{F}_{i,j_k}(k=1,2,...,O(C_{sp}))$, where the first subscript denotes the row cluster $i$. Different from the fill-ins at leaf level, the $\mathbf{F}_{i,j_k}$ now is of $(2k_{l+1} \times 2k_{l+1})$ size at a non-leaf level $l$. We then compute

$$\mathbf{G}_i = (\mathbf{I} - \mathbf{T}_i \mathbf{T}_i^H)(\sum_k \mathbf{F}_{i,j_k} \mathbf{F}_{i,j_k}^H)(\mathbf{I} - \mathbf{T}_i \mathbf{T}_i^H)^H, \quad (25)$$

where all the $k$ admissible blocks formed by cluster $i$ are considered. Performing an SVD on $\mathbf{G}_i$, and truncating the singular vector matrix based on prescribed accuracy $\epsilon_{fill-in}$, we obtain

$$\mathbf{G}_i \stackrel{\epsilon_{fill-in}}{=} (\mathbf{T}_i^{add}) \Sigma (\mathbf{T}_i^{add})^H. \quad (26)$$

The $\mathbf{T}_i^{add}$ is hence the additional column space we should supplement in the transfer matrix. The new transfer matrix $\widetilde{\mathbf{T}}_i$ can then be obtained as

$$\widetilde{\mathbf{T}}_i = [\mathbf{T}_i \ \mathbf{T}_i^{add}]. \quad (27)$$

The cost of this step for each cluster at level $l$ is $O(k_l^3)$.

**Step 1 at a non-leaf level.** We compute the complementary cluster basis $\widetilde{\mathbf{T}}_i^\perp$, which is orthogonal to $\widetilde{\mathbf{T}}_i$. This computation is not expensive, as it costs only $O(k_l^3)$ at level $l$. We then combine $\widetilde{\mathbf{T}}_i$ with $\widetilde{\mathbf{T}}_i^\perp$ to form $\widetilde{\mathbf{Q}}_i$ matrix as the following

$$\widetilde{\mathbf{Q}}_i = [(\widetilde{\mathbf{T}}_i^\perp)_{2k_{l+1} \times (2k_{l+1}-k_l)} \ (\widetilde{\mathbf{T}}_i)_{2k_{l+1} \times k_l}]. \quad (28)$$

**Step 2 at a non-leaf level.** With $\widetilde{\mathbf{Q}}_i$ computed, we build a block diagonal matrix

$$\mathbf{Q}_i = diag\{\mathbf{I}_1, \mathbf{I}_2, ..., \widetilde{\mathbf{Q}}_i, ..., \mathbf{I}_{\#\{\text{clusters at level } l\}}\}, \quad (29)$$

using which we construct

$$\mathbf{Z}_{cur} = \mathbf{Q}_i^H \mathbf{Z}_{cur} \overline{\mathbf{Q}}_i, \quad (30)$$

where $\mathbf{Z}_{cur}$ only consists of the block that remains to be factorized, as shown by the bottom right matrix block in Fig. 2(i), which is not an identity. Similar to that in the leaf level, by doing so, we zero out first $(2k_{l+1} - k_l)$ rows of the admissible blocks formed by row cluster $i$, as well as their transpose entries in $\mathbf{Z}_{cur}$, as shown by the white blocks in Fig. 2(j). The real computation of (30) is hence in the inadmissible blocks of non-leaf cluster $i$, the cost of which is $O(k_l)^3$.

**Step 3 at a non-leaf level.** We perform a partial LU factorization to eliminate the first $(2k_{l+1} - k_l)$ unknowns of cluster $i$ in $\mathbf{Z}_{cur}$. Let this set of unknowns be denoted by $i'$. We obtain

$$\mathbf{Z}_{cur} = \begin{bmatrix} \mathbf{L}_{i'i'} & \mathbf{0} \\ \mathbf{Z}_{ci'} \mathbf{U}_{i'i'}^{-1} & \mathbf{I} \end{bmatrix} \begin{bmatrix} \mathbf{I} & \mathbf{0} \\ \mathbf{0} & \widetilde{\mathbf{Z}}_{cc} \end{bmatrix} \begin{bmatrix} \mathbf{U}_{i'i'} & \mathbf{L}_{i'i'}^{-1} \mathbf{Z}_{i'c} \\ \mathbf{0} & \mathbf{I} \end{bmatrix} \quad (31)$$

where $\mathbf{F}_{i'i'} = \mathbf{L}_{i'i'} \mathbf{U}_{i'i'}$, and $\widetilde{\mathbf{Z}}_{cc}$ is a Schur complement. Each fill-in block is again either directly added, if the to-be added block is inadmissible, or stored temporarily in the corresponding admissible block. Eqn. (31) can further be written in short as

$$\mathbf{Z}_{cur} = \mathbf{L}_i^l \begin{bmatrix} \mathbf{I} & \mathbf{0} \\ \mathbf{0} & \widetilde{\mathbf{Z}}_{cc} \end{bmatrix} \mathbf{U}_i^l. \quad (32)$$

At the non-leaf level, each of $\mathbf{L}_{i'i'}$ and $\mathbf{U}_{i'i'}$ is of size at most $(2k_{l+1} - k_l)$ by $2k_{l+1}$. The $\mathbf{Z}_{i'c}$ ($\mathbf{Z}_{ci'}$) has only $C_{sp}$ blocks which are also inadmissible, each of which is of size at most $2k_{l+1}$ by $2k_{l+1}$. The schur complement involves $C_{sp}^2$ fill-in blocks, each of which is of size at most $2k_{l+1}$ by $2k_{l+1}$. Therefore, generating the L and U factors, as well as computing the fill-in blocks, only involves a constant number of computations, each of which scales as $O(k_l^3)$ at a non-leaf level. The storage of the L and U factors, as well as the fill-in blocks, only takes $O(k_l^2)$ multiplied by constant $O(C_{sp}^2)$.

*3) Updating the coupling matrix at current level and the transfer matrix at one level higher:* After the computation at the current non-leaf level $l$ is finished, before proceeding to the next level, $l-1$, we need to update the coupling matrices at current level, as well as the transfer matrices at one level higher. This is because the transfer matrices have been updated at the current level.

Similar to the procedure at the leaf level, for admissible blocks without fill-ins, their coupling matrices can be readily updated by appending zeros, as shown in (20). For admissible blocks with fill-ins, their coupling matrices are updated by adding the fill-in part onto the augmented coupling matrix

$\widetilde{\mathbf{S}}_{i,k}$ as
$$\widetilde{\mathbf{S}}_{i,k} = \widetilde{\mathbf{S}}_{i,k} + \widetilde{\mathbf{T}}_i^H \mathbf{F}_{i,k} \overline{\widetilde{\mathbf{T}}}_k, \tag{33}$$

in which $\mathbf{F}_{i,k}$ represents the fill-in block associated with the admissible block formed between clusters $i$ and $k$.

We also need to update the transfer matrices at one level higher, i.e., at the $l-1$ level. This can be readily done by adding zeros to the original transfer matrices, corresponding to the matrix block of $\mathbf{T}_i^{add}$.

### C. Pseudo-codes and Summary

The aforementioned computation at the non-leaf level is continued level by level until we reach the minimal level that has admissible blocks, denoted by $l_0$. The overall procedure can be realized by the pseudo-code shown in **Algorithm 1**.

---
**Algorithm 1** Proposed direct solution of general $\mathcal{H}^2$-matrices
---
1: Let $\mathbf{Z}_{cur} = \mathbf{Z}$
2: **for** $l = L$ to $l_0$ **do**
3:     **for** $cluster: i = 1$ to $2^l$ **do**
4:        Update cluster basis (Step 0)
5:        Obtain projection matrix $\mathbf{Q}_i$ (Step 1)
6:        Do $\mathbf{Z}_{cur} = \mathbf{Q}_i^H \mathbf{Z}_{cur} \overline{\mathbf{Q}}_i$ (Step 2)
7:        Perform partial LU $\mathbf{Z}_{cur} = \mathbf{L}_i^l \mathbf{Z}_{cnt} \mathbf{U}_i^l$ (Step 3)
8:        Let $\mathbf{Z}_{cur} = \mathbf{Z}_{cnt}$
9:     **end for**
10:    Update coupling matrices at level $l$
11:    Update transfer matrices at level $l-1$
12:    Permute and merge small $k_l \times k_l$ matrices to next level
13: **end for**
14: Do LU on the final block, $\mathbf{Z}_{cur} = \mathbf{LU}$.
---

The pseudo-code for Step 0 and Step 1 is shown in Algorithm 2. For each cluster, we check whether its associated admissible blocks have fill-ins or not. If it has fill-ins, we first update the cluster basis based on (17) or transfer matrix based on (26), we then obtain (19) or (28).

---
**Algorithm 2** Steps 0 and 1: Update $\mathbf{V}_i$ and Obtain $\mathbf{Q}_i$
---
1: **if** cluster $i$ has fill-ins **then**
2:     Update $\mathbf{V}_i$ or $\mathbf{T}_i$
3: **end if**
4: Perform QL on $\mathbf{V}_i$ ($\mathbf{T}_i$) to find $\mathbf{V}_i^\perp$ ($\mathbf{T}_i^\perp$).
5: Combine $\mathbf{V}_i$ ($\mathbf{T}_i$) and $\mathbf{V}_i^\perp$ ($\mathbf{T}_i^\perp$) to get $\mathbf{Q}_i$.
---

After getting $\mathbf{Q}_i$ matrix, we use it to zero out entries in the admissible blocks related to cluster $i$, which is Step 2. At this step, the operation of (5) is essentially to left multiply the $i$-th block row of $\mathbf{Z}_{cur}$ using $\widetilde{\mathbf{Q}}_i^H$, and right multiply the $i$-th block column of $\mathbf{Z}_{cur}$ using $\overline{\widetilde{\mathbf{Q}}}_i$. If the block being multiplied is an inadmissible block, we perform a full matrix multiplication; if the block is an admissible block, no computation is needed as the product is known to be the coupling matrix of this block, which is either the original one or the one superposed with the fill-in's contribution shown in (22) and (33). The pseudo-code of Step 2 is shown in Algorithm 3.

---
**Algorithm 3** Step 2: Do $\mathbf{Z}_{cur} = \mathbf{Q}_i^H \mathbf{Z}_{cur} \overline{\mathbf{Q}}_i$
---
1: **for** Column cluster $j$ forming a block with row $i$ **do**
2:     **if** $\mathbf{Z}_{cur,ij}$ is inadmissible **then**
3:        Do $\widetilde{\mathbf{Q}}_i^H \mathbf{Z}_{cur,ij}$.
4:     **else**
5:        No computation needed.
6:     **end if**
7: **end for**
8: **for** Row cluster $j$ forming a block with column $i$ **do**
9:     **if** $\mathbf{Z}_{ji}$ is inadmissible **then**
10:       Do $\mathbf{Z}_{cur,ji} \overline{\widetilde{\mathbf{Q}}}_i$.
11:    **else**
12:       No computation needed.
13:    **end if**
14: **end for**
---

Step 3 is shown in Algorithm 4, where we perform a partial LU factorization on diagonal matrix block $\mathbf{F}_{i'i'}$. Step 4 in this algorithm is realized by generating $O(C_{sp}^2)$ fill-in blocks and adding (or storing) them upon the corresponding blocks in current matrix being factorized, i.e., $\mathbf{Z}_{cur}$.

---
**Algorithm 4** Step 3: Perform $\mathbf{Z}_{cur} = \mathbf{L}_i^l \mathbf{Z}_{cnt} \mathbf{U}_i^l$
---
1: Partial LU on diagonal matrix block $\mathbf{Z}_{i'i'} = \mathbf{L}_{i'i'} \mathbf{U}_{i'i'}$.
2: Update $\mathbf{U}_{i'c} \leftarrow \mathbf{L}_{i'i'}^{-1} \mathbf{Z}_{i'c}$
3: Update $\mathbf{L}_{ci'} \leftarrow \mathbf{Z}_{ci'} \mathbf{U}_{i'i'}^{-1}$
4: Update $\mathbf{Z}_{cnt} \leftarrow \mathbf{Z}_{cc} - \mathbf{Z}_{ci'} \mathbf{U}_{i'i'}^{-1} \mathbf{L}_{i'i'}^{-1} \mathbf{Z}_{i'c}$
---

The aforementioned direct solution procedure results in the following factorization of $\mathbf{Z}$:
$$\begin{aligned} \mathbf{Z} &= \mathcal{L}\mathcal{U}, \\ \mathcal{L} &= \left\{ \prod_{l=L}^{l_0} \left[ \left( \prod_{i=1}^{2^l} \mathbf{Q}_i^l \mathbf{L}_i^l \right) \mathbf{P}_l^T \right] \right\} \mathbf{L}, \\ \mathcal{U} &= \mathbf{U} \left\{ \prod_{l=l_0}^{L} \left[ \mathbf{P}_l \left( \prod_{i=2^l}^{1} \mathbf{U}_i^l \mathbf{Q}_i^{l^T} \right) \right] \right\}. \end{aligned} \tag{34}$$

In the above, $\mathbf{Q}_i^l$ denotes $\mathbf{Q}_i$ at level $l$. For leaf level, it is shown in (6), and for a non-leaf level, it is shown in (29). Each $\mathbf{Q}_i^l$ has only one block that is not identity, whose size is $2k_{l+1} \times 2k_{l+1}$ at the non-leaf level, and $leafsize \times leafsize$ at the leaf level. Each $\mathbf{L}_i^l$ is shown in (13) or (32), which has $O(C_{sp})$ blocks of $leafsize$ at the leaf level, and $O(C_{sp})$ blocks of $O(k_{l+1})$ size at a non-leaf level. The same is true to $\mathbf{U}_i^l$. The $\mathbf{L}$ and $\mathbf{U}$ without sub- and super-scripts are the $\mathbf{L}$ and $\mathbf{U}$ factors of the final matrix that remains to be factorized in level $l_0$. The non-identity blocks in $\mathbf{L}$ and $\mathbf{U}$ are of size $O(2^{l_0} k_{l_0})$. The $\mathbf{P}_l$ matrix is a permutation matrix that merges small $k_l \times k_l$ matrices at level $l$ to one level higher.

The factorization shown in (34) also results in an explicit

(a) $\mathbf{L}_1$ matrix.  (b) $\mathbf{U}_1$ matrix.

Fig. 3: (a) $\mathbf{L}_1$ matrix. (b) $\mathbf{U}_1$ matrix.

form of $\mathbf{Z}$'s inverse, which can be written as

$$\mathbf{Z}^{-1} = \left\{\prod_{l=L}^{l_0}\left[\left(\prod_{i=1}^{2^l}\overline{\mathbf{Q}}_i^l(\mathbf{U}_i^l)^{-1}\right)\mathbf{P}_l^T\right]\right\}\mathbf{U}^{-1}\times \quad (35)$$
$$\mathbf{L}^{-1}\left\{\prod_{l=l_0}^{L}\left[\mathbf{P}_l\left(\prod_{i=2^l}^{1}(\mathbf{L}_i^l)^{-1}\mathbf{Q}_i^{l\,H}\right)\right]\right\}.$$

In addition, in the proposed direct solution, one can also flexibly stop at a level before $l_0$ is reached. This may be desired if at a tree level, the number of zeros introduced is small as compared to the block size at that level. Although one can still proceed to the next level, the efficiency gain may be smaller than that gained if one just stops at the current level, and finishes the factorization. This typically happens at a few levels lower than the $l_0$ level.

## IV. PROPOSED MATRIX SOLUTION (BACKWARD AND FORWARD SUBSTITUTION) ALGORITHM

After factorization, the solution of $\mathbf{Z}x = b$ can be obtained in two steps: $\mathcal{L}y = b$ which can be called as a backward substitution, and $\mathcal{U}x = y$ called as a forward substitution. Rigorously speaking, since projection matrix $\mathbf{Q}_i^l$ and permutation matrix $\mathbf{P}_l$ are involved in $\mathcal{L}$ and $\mathcal{U}$ factors, the matrix solution here is not a conventional forward/backward substitution.

### A. Forward Substitution

In this step, we solve

$$\mathcal{L}y = b. \quad (36)$$

From (35), it can be seen that the above can be solved as

$$y = \mathcal{L}^{-1}b = \mathbf{L}^{-1}\left\{\prod_{l=l_0}^{L}\left[\mathbf{P}_l\left(\prod_{i=2^l}^{1}(\mathbf{L}_i^l)^{-1}\mathbf{Q}_i^{l\,H}\right)\right]\right\}b. \quad (37)$$

Take the 3-level $\mathcal{H}^2$-matrix as an example, the above is nothing but

$$y = \mathbf{L}^{-1}\mathbf{P}_1\mathbf{L}_4^{-1}\mathbf{Q}_4^H\mathbf{L}_3^{-1}\mathbf{Q}_3^H\mathbf{L}_2^{-1}\mathbf{Q}_2^H\mathbf{L}_1^{-1}\mathbf{Q}_1^H b. \quad (38)$$

Obviously, three basic operations are involved in the forward substitution: multiplying $(\mathbf{Q}_i^l)^H$ by a vector, multiplying $(\mathbf{L}_i^l)^{-1}$ by a vector, and $\mathbf{P}_l$-based permutation.

The algorithm of the forward substitution is shown in **Algorithm 5**. The $(\mathbf{L}_i^l)^{-1}b$ in step 4 of this algorithm only involves $O(C_{sp})$ matrix-vector multiplications of $leafsize$ or rank size associated with cluster $i$ at level $l$. To see this point more clearly, we can use the $\mathbf{L}_1$ matrix shown in Fig. 3(a) as an example, where $\mathbf{L}_1$ is resultant from the partial factorization of cluster 1 in the $\mathcal{H}^2$-matrix shown in Fig. 1(b). It can be seen that $\mathbf{L}_1$ has only $O(C_{sp})$ blocks that need to be computed and stored, each of which is $leafsize$ at the leaf level, or rank size at a non-leaf level $l$. These blocks are the inadmissible blocks formed with cluster $i$. For such type of $\mathbf{L}_i$, the $\mathbf{L}_i^{-1}b$ can be readily computed. For $\mathbf{L}_1$ shown in Fig. 3(a), the $\mathbf{L}_1^{-1}b$ is

$$\mathbf{L}_1^{-1}b = \begin{bmatrix}\begin{bmatrix}\mathbf{L}_{11}^{-1} & \mathbf{0}\\ -\mathbf{T}_{11}\mathbf{L}_{11}^{-1} & \mathbf{I}\end{bmatrix} & \mathbf{0} & \mathbf{0} & \mathbf{0}\\ \begin{bmatrix}-\mathbf{T}_{21}\mathbf{L}_{11}^{-1} & \mathbf{0}\end{bmatrix} & \mathbf{I} & \mathbf{0} & \mathbf{0}\\ \begin{bmatrix}-\mathbf{T}_{31}\mathbf{L}_{11}^{-1} & \mathbf{0}\end{bmatrix} & \mathbf{0} & \mathbf{I} & \mathbf{0}\\ \mathbf{0} & \mathbf{0} & \mathbf{0} & \mathbf{I}\end{bmatrix}\begin{bmatrix}b_{11}\\ b_{12}\\ b_2\\ b_3\\ b_4\end{bmatrix}$$
$$= \begin{bmatrix}\widetilde{b}_{1,1}\\ b_{1,2} - \mathbf{T}_{11}\widetilde{b}_{1,1}\\ b_2 - \mathbf{T}_{21}\widetilde{b}_{1,1}\\ b_3 - \mathbf{T}_{31}\widetilde{b}_{1,1}\\ b_4\end{bmatrix}$$
(39)

where $\mathbf{T}_{i1} = \mathbf{F}_{i1}\mathbf{U}_{11}^{-1}$ ($i = 1, 2, 3$), and $\widetilde{b}_{1,1} = \mathbf{L}_{11}^{-1}b_{1,1}$. Therefore, this operation actually only involves $O(C_{sp})$ matrix-vector multiplications related to clusters 1, 2 and 3, which are those clusters that form inadmissible blocks with cluster 1, the number of which is bounded by $C_{sp}$. After forward substitution, the right hand side $b$ vector is overwritten by solution $y$.

**Algorithm 5** Forward Substitution $b \leftarrow \mathcal{L}^{-1}b$
1: **for** $l = L$ to $l_0$ **do**
2:    **for** $cluster : i = 1$ to $2^l$ **do**
3:       Update $b$'s $b_i$ part by $b_i \leftarrow (\widetilde{\mathbf{Q}}_i^l)^H b_i$.
4:       Do $b \leftarrow (\mathbf{L}_i^l)^{-1}b$
5:    **end for**
6:    Permute $b$ by multiplying with $\mathbf{P}_l$
7: **end for**
8: Do $b \leftarrow \mathbf{L}^{-1}b$

### B. Backward Substitution

After finishing forward substitution, we solve

$$\mathcal{U}x = y \quad (40)$$

From (35), it can be seen that the above can be solved as

$$x = \mathcal{U}^{-1}y = \left\{\prod_{l=L}^{l_0}\left[\left(\prod_{i=1}^{2^l}\overline{\mathbf{Q}}_i^l{\mathbf{U}_i^l}^{-1}\right)\mathbf{P}_l^T\right]\right\}\mathbf{U}^{-1}y. \quad (41)$$

For the same 3-level $\mathcal{H}^2$-matrix example, the above is nothing but

$$x = \overline{\mathbf{Q}}_1\mathbf{U}_1^{-1}\overline{\mathbf{Q}}_2\mathbf{U}_2^{-1}\overline{\mathbf{Q}}_3\mathbf{U}_3^{-1}\overline{\mathbf{Q}}_4\mathbf{U}_4^{-1}\mathbf{P}_1^T\mathbf{U}^{-1}y. \quad (42)$$

The backward substitution starts from minimal admissible block level $l_0$, and finishes at the leaf level, as shown in **Algorithm 6**. The $(\mathbf{U}_i^l)^{-1}y$ in step 5 is performed efficiently similar to how (39) is used to compute $(\mathbf{L}_i^l)^{-1}b$. There are also three basic operations in **Algorithm 6**: $(\mathbf{U}_i^l)^{-1}$-based matrix-vector multiplication, $\overline{\mathbf{Q}}_i^l$-based matrix-vector multiplication, and $\mathbf{P}_l$-based permutation. After the backward substitution is done, $y$ is overwritten by solution $x$.

---

**Algorithm 6** Backward Substitution $y \leftarrow \mathcal{U}^{-1}y$

1: Do $y \leftarrow \mathbf{U}^{-1}y$
2: **for** $l = l_0$ **to** $L$ **do**
3:    Permute $y$ by front multiplying $\mathbf{P}_l^T$
4:    **for** $cluster i = 2^l$ **to** $1$ **do**
5:      Do $y \leftarrow (\mathbf{U}_i^l)^{-1}y$
6:      Update $y_i \leftarrow (\overline{\mathbf{Q}}_i^l)y_i$.
7:    **end for**
8: **end for**

---

## V. ACCURACY AND COMPLEXITY ANALYSIS

In this section, we analyze the accuracy and computational complexity of the proposed direct solver.

### A. Accuracy

Different from [9], [11], [15], in the proposed new direct solvers, the accuracy is directly controlled. When generating an $\mathcal{H}^2$-matrix to represent the original dense matrix, the accuracy is controlled by $\epsilon_{\mathcal{H}^2}$, which is the same as in [9], [11], [15]. The accuracy of the proposed direct solution is controlled by $\epsilon_{fill-in}$. As can be seen, there are no approximations involved in the proposed solution procedure, except that when we update the cluster bases to account for the contribution from the fill-ins, we use (17) or (26) to truncate the singular vector matrix. The accuracy of this step is controlled by accuracy parameter $\epsilon_{fill-in}$.

### B. Time and Memory Complexity

From (34), it is evident that the whole computation involves $\mathbf{Q}_i^l$, $\mathbf{L}_i^l$, and $\mathbf{U}_i^l$ for each cluster $i$ at every level $l$, the storage and computation of each of which are respectively $O(k_l^2)$, and $O(k_l^3)$ at a non-leaf level, and constant at the leaf level. Recall that $\mathbf{Q}_i^l$ is obtained from Step 0 and Step 1, $\mathbf{L}_i^l$, and $\mathbf{U}_i^l$ are obtained from Step 3, whose underlying blocks are generated in Step 2. As for the permutation matrix $\mathbf{P}_l$, it is sparse involving only $O(2^l)$ integers to store at level $l$, and hence having an $O(N)$ cost in both memory and time.

From another perspective, the overall procedure of the proposed direct solution is shown in **Algorithm 1**. It involves $O(L)$ levels of computation. At each level, there are $2^l$ clusters. For each cluster, we have performed four steps of computation from Step 0, 1, 2, to 3, as shown by the tasks listed at the end of lines 4–7 in **Algorithm 1**. Each of these steps costs $O(k_l)^3$ for each cluster, as analyzed in Section III. After each level is done, we update coupling matrices, and transfer matrices at one-level higher, the cost of which is also $O(k_l)^3$ for every admissible block, or cluster.

Taking into account the rank's growth with electrical size [21], the factorization and inversion complexity of the proposed direct solutions can be found for electrically large volume IEs as

$$\sum_{l=0}^{L} C_{sp}^2 2^l k_l^3 = C_{sp}^2 \sum_{l=0}^{L} 2^l \left[\left(\frac{N}{2^l}\right)^{\frac{1}{3}}\right]^3 = O(NlogN), \quad (43)$$

since at every level, there are $2^l$ clusters, and each cluster is associated with $O(C_{sp}^2)$ computations, and each of which costs $O(k_l^3)$. Here, $k_l$ is proportional to the cubic root of block size at level $l$, because it is the electric size at level $l$ in an electrically large volume IE. For small or modest electric sizes, the rank can be well approximated by a constant, and hence the resulting time complexity is strictly $O(N)$. The solution time and memory only involve $O(k_l^2)$ computations for each cluster. Hence, their complexities are

$$\sum_{l=0}^{L} C_{sp} 2^l k_l^2 = C_{sp} \sum_{l=0}^{L} 2^l \left[\left(\frac{N}{2^l}\right)^{\frac{1}{3}}\right]^2 = O(N), \quad (44)$$

which is linear irrespective of electric size.

## VI. NUMERICAL RESULTS

In order to demonstrate the accuracy and low computational complexity of the proposed direct solution of general $\mathcal{H}^2$-matrices, we apply it to solve electrically large volume integral equations (VIE) for electromagnetic analysis. The VIE formulation we use is based on [22] with SWG vector bases for expanding electric flux density in each tetrahedral element. First, the scattering from a dielectric sphere, whose analytical solution is known, is simulated to validate the proposed direct solver. Then, a variety of large-scale examples involving over one million unknowns are simulated on a single CPU core to examine the accuracy and complexity of the proposed direct solver. The $\mathcal{H}^2$-matrix for each example is constructed based on the method described in [14], [15]. The accuracy parameter $\epsilon_{\mathcal{H}^2}$ used is $10^{-3}$ for the $\mathcal{H}^2$-matrix construction. The accuracy parameter $\epsilon_{fill-in}$ used in the direct solution varies in different examples. The computer used has an Intel(R) Xeon(R) CPU E5-2690 v2 running at 3.00GHz, and only a *single core* is employed for carrying out the computation.

### A. Scattering from a Dielectric Sphere

We simulate a dielectric sphere of radius $a$ with different electric sizes $k_0 a$ and dielectric constant $\epsilon_r$. The incident electric field is a normalized $-z$ propagating plane wave polarized along the $x-$direction. This incident field is also used in the following examples. During the $\mathcal{H}^2$-matrix construction, $leafsize = 25$ and $\eta = 1$ are used. The accuracy parameter in the direct solution is set to be $\epsilon_{fill-in} = 1e - 6$. In Fig. 4, the radar-cross-section (RCS in dBsm) of the sphere is plotted as a function of $\theta$ for zero azimuth for $k_0 a = 0.408, \epsilon_r = 36$; $k_0 a = 0.408, \epsilon_r = 4$; $k_0 a = 0.816, \epsilon_r = 4$; and $k_0 a = 1.632, \epsilon_r = 4$ respectively. They all show good agreement with the analytical Mie Series solution, which validates the proposed direct solution algorithm.

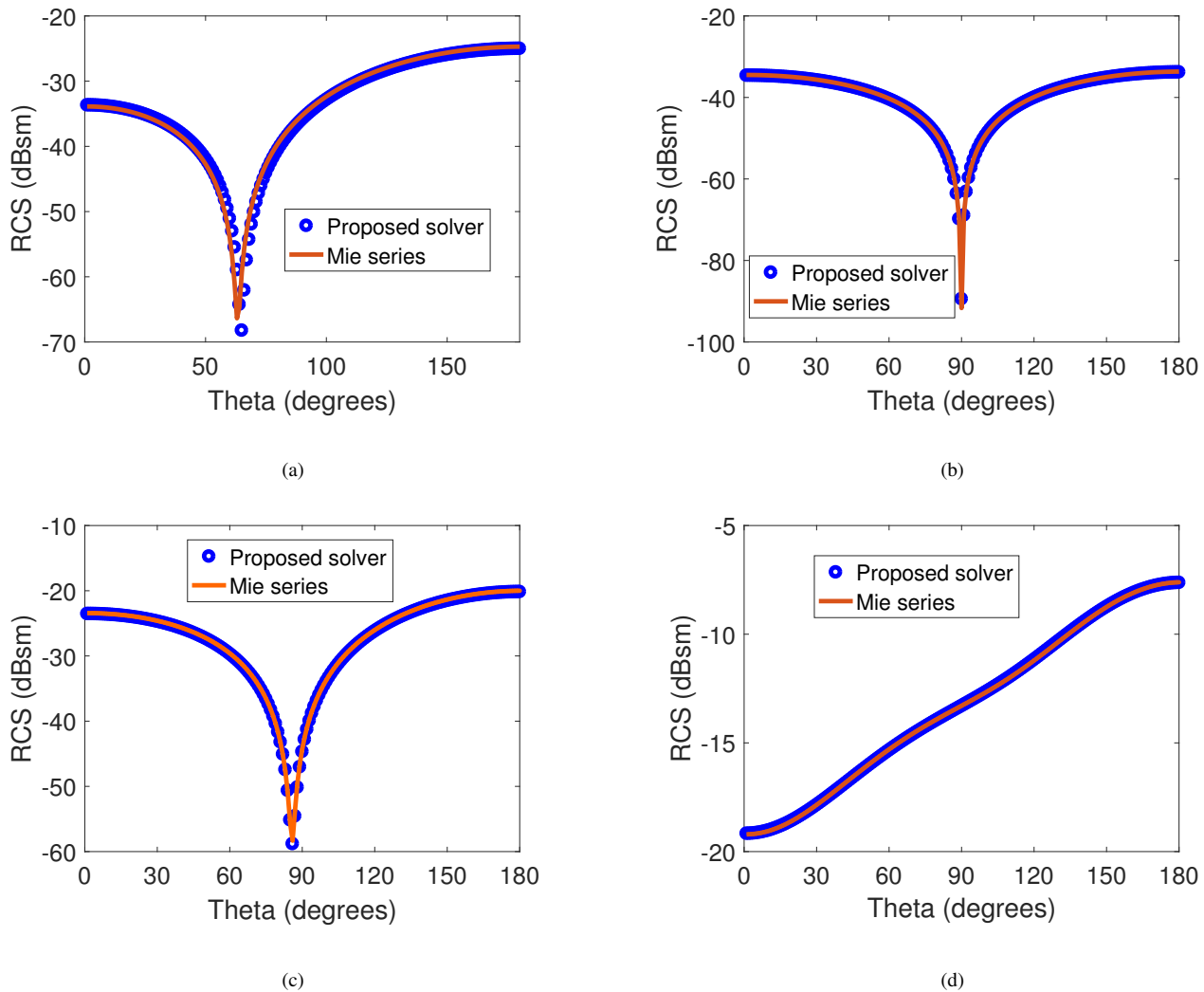

Fig. 4: Simulated RCS of a dielectric sphere for different electric sizes and dielectric constants. (a) $k_0 a = 0.408$, $\epsilon_r = 36.0$. (b) $k_0 a = 0.408$, $\epsilon_r = 4.0$. (c) $k_0 a = 0.816$, $\epsilon_r = 4.0$. (d) $k_0 a = 1.632$, $\epsilon_r = 4.0$.

### B. Large-scale Dielectric Rod

A dielectric rod with $\epsilon_r = 2.54$ is simulated at 300 MHz. The cross-sectional area of the rod is a square of side length $0.4\lambda$. The length of the rod is increased from $20\lambda$ to $485\lambda$ to study the complexity of the proposed direct solver versus $N$, and thereby electric size. The mesh size is fixed to be $0.1\lambda$ for all these rod examples. The increase of the electric size of the rod length results in $N$ increasing from 41,632, 166,432, 332,832, 665,632, to 1,008,832. The last case has over one million unknowns. The $leaf size = 25$ and $\eta = 1$ are used in the $\mathcal{H}^2$-matrix construction. The accuracy parameter $\epsilon_{fill-in}$ is set to be $10^{-4}$ in the proposed direct solution.

Since only one dimension of the simulated structure is changing, this example represents 1-D problem characteristics [21], whose rank remains constant irrespective of the electric size. In Fig. 5a, we plot the memory cost for storing the original $\mathcal{H}^2$-matrix, as well as the memory cost of the matrix factorization, as a function of $N$. These two are different, since there are additional matrices generated during the factorization procedure such as the $\mathbf{Q}_i^l$, the $\mathbf{P}_i^l$ matrices, etc. Obviously, both the matrix storage and the factorization memory scale linearly with the number of unknowns, which verifies our theoretical analysis given in (44). The factorization time and solution time with respect to $N$ are plotted in Fig. 5b. They are also shown to scale linearly with $N$ instead of being $O(N log N)$, because of the constant rank in this 1-D example. The error of the proposed direct solution is measured by computing the relative residual $\epsilon_{rel} = ||\mathbf{Z}_{\mathcal{H}^2} x - b||/||b||$, where $\mathbf{Z}_{\mathcal{H}^2}$ is the input $\mathcal{H}^2$-matrix to be solved. This error is listed in Table I. Excellent accuracy can be observed in the entire unknown range.

TABLE I: Direct solution error measured by relative residual for the dielectric rod example.

| $N$ | 4,1632 | 166,432 | 332,832 | 665,632 | 1,008,832 |
|---|---|---|---|---|---|
| $\epsilon_{rel}$ | 0.57% | 0.68% | 0.48% | 0.47% | 0.30% |

### C. Large-scale Dielectric Slab

We then simulate a dielectric slab with $\epsilon_r = 2.54$ at 300 MHz, which is also simulated in [14], [15]. The thickness of the slab is fixed to be $0.1\lambda_0$. The width and length are simultaneously increased from $4\lambda_0$, $8\lambda_0$, $16\lambda_0$, to $32\lambda_0$. With a mesh size of $0.1\lambda_0$, the resultant $N$ ranges from 22,560 to 1,434,880 for this suite of slab structures. The $leafsize$ is chosen to be 25 and $\eta = 1$. The $\epsilon_{fill-in}$ is set to be $10^{-5}$.

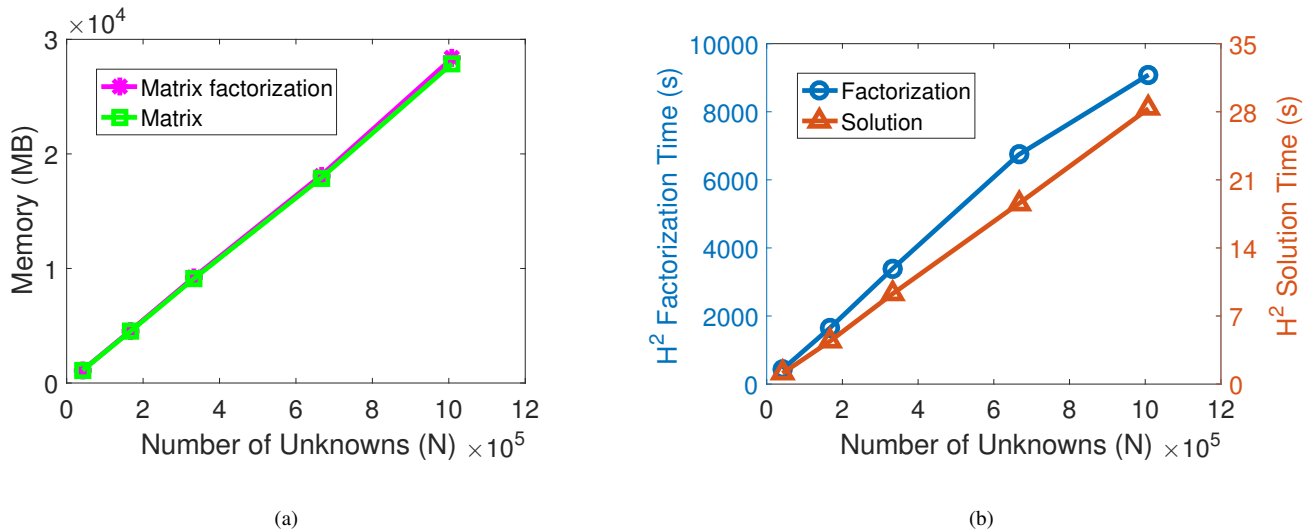

Fig. 5: Solver performance of a suite of dielectric rod examples. (a) Memory scaling v.s. $N$. (b) Factorization time and solution time scaling v.s. $N$.

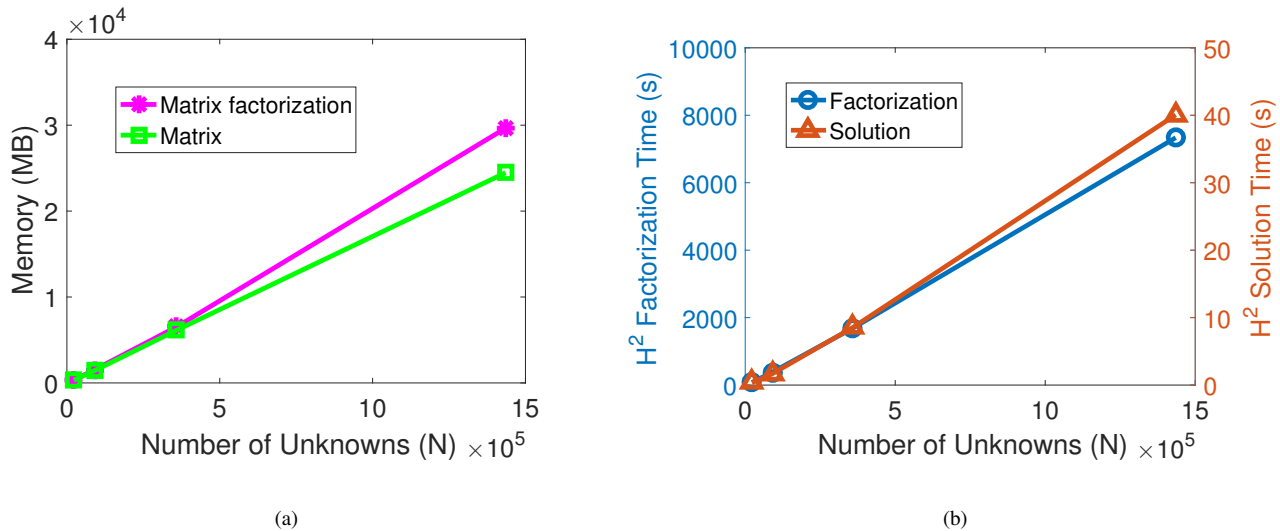

Fig. 6: Solver performance of a large-scale dielectric slab. (a) Memory scaling v.s. $N$. (b) Factorization time and solution time scaling v.s. $N$.

In Fig. 6a, we plot the memory cost of the proposed direct solver as a function of $N$. In Fig. 6b, we plot the factorization time, and the solution time as a function of $N$. Clear linear complexity is observed in both memory consumption and CPU run time. Notice that this example represents a 2-D problem characteristics. Based on [21], the rank's growth rate with electric size is lower than linear, and being a square-root of the log-linear of the electric size. Substituting such a rank's growth rate into the complexity analysis shown in (44) and (43), we will obtain linear complexity in both memory and time. The relative residual $\epsilon_{rel}$ of the proposed direct solution is listed in Table II. Again, excellent accuracy can be observed in the entire unknown range.

TABLE II: Direct solution error measured by relative residual for the dielectric slab example.

| $N$ | 22,560 | 89,920 | 359,040 | 1,434,880 |
|---|---|---|---|---|
| $\epsilon_{rel}$ | 0.085% | 0.15% | 0.58% | 0.37% |

Since this example is also simulated in our previous work [14], [15], we have compared the two direct solvers in CPU run time and accuracy. As can be seen from Table III, the proposed new direct solution takes much less time than that of [14], [15]. The speedup is more than one order of magnitude in large examples. Even though the factorization time is shown, as can be seen from (35), the inversion time in the proposed algorithm is similar to factorization time. In addition, because of a direct error control, the error of the proposed solution shown in Table II is also much less than that of [14], [15], which is shown to be 3.9%, and 7.49% for $N = 89,920$, and 359,040 respectively.

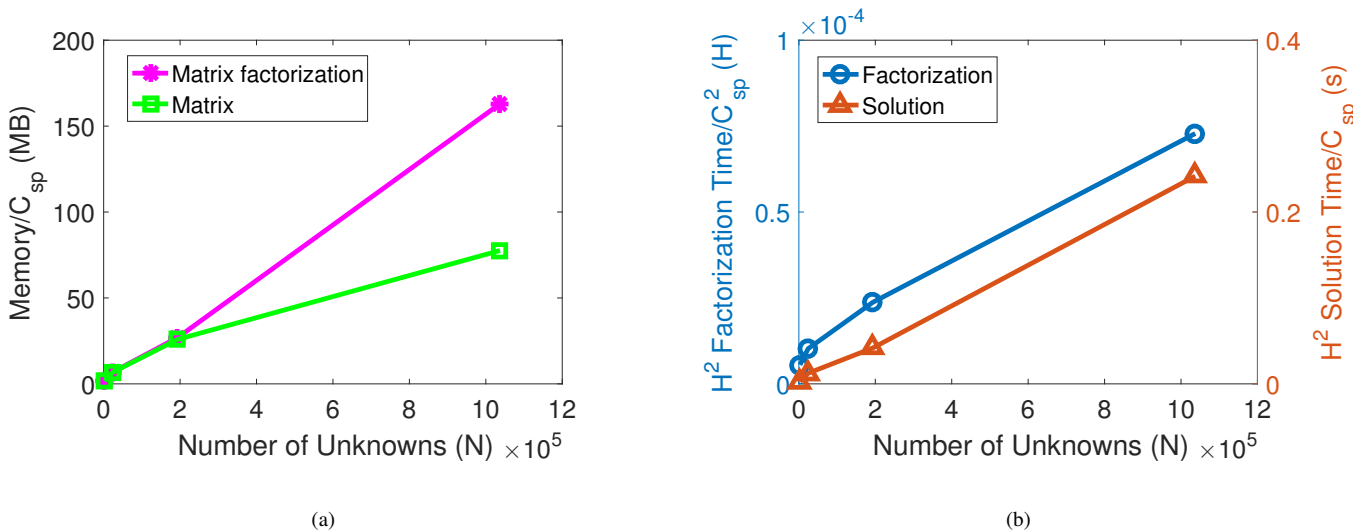

Fig. 7: Solver performance of a suite of cube array structures. (a) Memory scaling v.s. $N$. (b) Factorization time and solution time scaling v.s. $N$.

TABLE III: Performance comparison with [14], [15] for the dielectric slab example.

| $N$ | 89,920 | 359,040 | 1,434,880 |
|---|---|---|---|
| Factorization (s) [This] | 380 | 1,690 | 7,335 |
| Solution (s) [This] | 1.66 | 8.55 | 40.06 |
| Inversion (s) [15] | 2.75e+03 | 1.65e+04 | 9.31e+04 |

### D. Large-scale Array of Dielectric Cubes

Next, we simulate a large-scale array of dielectric cubes at 300 MHz. The relative permittivity of the cube is $\epsilon_r = 4.0$. Each cube is of size $0.3\lambda_0 \times 0.3\lambda_0 \times 0.3\lambda_0$. The distance between adjacent cubes is kept to be $0.3\lambda_0$. The number of the cubes is increased along the $x$-, $y$-, and $z$-directions simultaneously from 2 to 14, thus producing a 3-D cube array from $2 \times 2 \times 2$ to $14 \times 14 \times 14$ elements. The number of unknowns $N$ is respectively 3,024, 24,192, 193,536, and 1,037,232 for these arrays. The $leafsize$ is 25 and $\eta = 1$. The $\epsilon_{fill-in}$ used is $10^{-5}$.

For the cubic growth of unknowns for 3-D problems, we observe that constant $C_{sp}$ is not saturated until in the order of millions of unknowns, as can be seen from Table IV. It is thus important to analyze the performances of the proposed direct solver as (Memory/Solution time)/$C_{sp}$ and (Factorization time)/$C_{sp}^2$ respectively. In Fig. 7a and Fig. 7b, we plot the storage cost normalized with $C_{sp}$, and the direct factorization/solution time divided by $C_{sp}^2/C_{sp}$ with respect to $N$. As can be seen, their scaling rate with $N$ agrees very well with our theoretical complexity analysis. The relative residual

TABLE IV: $C_{sp}$ as a function of $N$ for the dielectric cube array.

| $N$ | 3,024 | 24,192 | 193,536 | 378,000 | 1,037,232 |
|---|---|---|---|---|---|
| $C_{sp}$ | 16 | 42 | 95 | 327 | 270 |

$\epsilon_{rel}$ of the proposed direct solution is listed in Table V for this example, which again reveals excellent accuracy. For the over one-million unknown case which is a $14 \times 14 \times 14$ cube array having thousands of cube elements, the error is still controlled to be as small as 1.3%. The matrix factorization of such a large-scale example is completed in 19,118 s (approximately 5 hours), and one solution time only takes 65.28 s.

TABLE V: Direct solution error measured by relative residual for the cube array example.

| $N$ | 3,024 | 24,192 | 193,536 | 1,037,232 |
|---|---|---|---|---|
| $\epsilon_{rel}$ | 0.09% | 0.15% | 0.36% | 1.3% |

## VII. Conclusions

In this paper, we develop fast factorization and inversion algorithms for general $\mathcal{H}^2$-matrices. The key feature of this new direct solution is its directly controlled accuracy. The entire direct solution does not involve any approximation, except that the fill-in block is compressed to update the cluster basis. However, this compression is well controlled by accuracy parameter $\epsilon_{fill-in}$. The resultant direct solution has been applied to solve electrically large VIEs whose kernel is oscillatory and complex-valued, which has been difficult to solve. Clear $O(N)$ complexity is achieved in memory and solution time, and $O(NlogN)$ complexity is achieved for direct factorization and inversion irrespective of electric size. Millions of unknowns are directly solved on a single core CPU in fast CPU run time. Comparisons with state-of-the-art $\mathcal{H}^2$-based direct VIE solvers have demonstrated the accuracy of this new direct solution as well as a significantly improved computational efficiency.